\documentclass[12pt]{article}
\usepackage[a4paper,margin=1.0in]{geometry}
\usepackage[colorlinks,citecolor=magenta,linkcolor=black]{hyperref}
\pdfpagewidth=\paperwidth \pdfpageheight=\paperheight
\usepackage{amsfonts,amssymb,amsthm,amsmath,eucal,tabu,url}
\usepackage{pgf}
 \usepackage{array}
 \usepackage{tikz-cd}
 \usepackage{pstricks}
 \usepackage{pstricks-add}
 \usepackage{pgf,tikz}
 \usetikzlibrary{automata}
 \usetikzlibrary{arrows}
 \usepackage{indentfirst}
\usepackage{tabularx} 

\theoremstyle{plain}
\newtheorem{thm}{Theorem}[section]
\newtheorem{theorem}[thm]{Theorem}
\newtheorem*{theoremA}{Theorem A}
\newtheorem*{theoremB}{Theorem B}

\theoremstyle{definition}
\newtheorem{definition}[thm]{Definition}
\newtheorem{remark}[thm]{Remark}
\newtheorem{example}[thm]{Example}

\newtheorem{thevarthm}[thm]{\varthmname}

\newenvironment{varthm*}[1]{\trivlist\item[]{\bf #1.}\it}{\endtrivlist}

\renewcommand\geq{\geqslant}

\renewcommand\leq{\leqslant}

\newcommand\be{\begin{eqnarray*}}
\newcommand\ee{\end{eqnarray*}}

\newcommand\newop[2]{\def#1{\mathop{\rm #2}\nolimits}}
\newop\log{log}
\newop\ord{ord}
\newop\Gal{Gal}
\newop\SL{SL}
\newop\Bl{Bl}
\newop\mult{mult}
\newop\mass{mass}
\newop\div{div}
\newop\codim{codim}
\newop\sing{sing}
\newop\vdim{vdim}
\newop\edim{edim}
\newop\Ass{Ass}
\newop\size{size}
\newop\reg{reg}
\newop\satdeg{satdeg}
\newop\supp{supp}
\newop\Neg{Neg}
\newop\Nef{Nef}
\newop\Nefh{Nef_H}
\newop\Eff{Eff}
\newop\Zar{Zar}
\newop\MB{MB}
\newop\MBxC{MB\mathit{(x,C)}}
\newop\NnB{NnB}
\newop\Bigg{Big}
\newop\Effbar{\overline{\Eff}}

\def\keywordname{{\bfseries Keywords}}%
\def\keywords#1{\par\addvspace\medskipamount{\rightskip=0pt plus1cm
\def\and{\ifhmode\unskip\nobreak\fi\ $\cdot$
}\noindent\keywordname\enspace\ignorespaces#1\par}}
\def\subclassname{{\bfseries Mathematics Subject Classification
(2020)}\enspace}
\def\subclass#1{\par\addvspace\medskipamount{\rightskip=0pt plus1cm
\def\and{\ifhmode\unskip\nobreak\fi\ $\cdot$
}\noindent\subclassname\ignorespaces#1\par}}

\begin{document}
\title{$\mathcal{Q}$-conic arrangements in the complex projective plane}
\author{Piotr Pokora}
\date{\today}
\maketitle
\thispagestyle{empty}
\begin{abstract}
We study the geometry of $\mathcal{Q}$-conic arrangements in the complex projective plane. These are arrangements consisting of smooth conics and they admit certain quasi-homogeneous singularities. We show that such $\mathcal{Q}$-conic arrangements are never free. Moreover, we provide combinatorial constraints of the weak combinatorics of such arrangements.
\subclass{14C20, 32S22, 14N20}
\end{abstract}
\section{Introduction}
In the last decades there is a constantly growing interest in hyperplane arrangements in the complex projective spaces that are free. The reason standing behind this phenomenon is the celebrated Terao's freeness conjecture which asserts that the freeness of hyperplane arrangements over an arbitrary fixed field $\mathbb{F}$ is purely determined by the combinatorics, i.e., the intersection poset of a given arrangement. This conjecture is widely open, even in the case of line arrangements in the projective plane over $\mathbb{F}$. In this context, it is worth recalling a recent result due to Barakat and K\"uhne \cite{BarKuh} which tells us that the Terao's freeness conjecture holds with up to $14$ lines, and this result is characteristic free. Even if we have several deep results devoted to the freeness of line arrangements, we still try to understand specific examples that can potentially lead to a counterexample to the Terao's freeness conjecture. On the other hand, the freeness problem for reduced curves that are not line arrangements is even less understood. We have several results devoted to arrangements of reduced curves, mostly arrangements of smooth conics and lines. In \cite{SchenckToh}, Schenck and Toh\v{a}neanu study conic-line arrangements with quasi-homogeneous singularities and they provide, among many interesting results, an addition-deletion theorem. On the contrast, the authors also present a counterexample to a naive generalization of the Terao's freeness conjecture in the class of conic-line arrangements. Recently in \cite{DimcaPokora}, the author with Dimca have provided a complete classification of free conic-line arrangements in the complex projective plane with nodes, tacnodes, and ordinary triple points -- it turns out that we have, up to the projective equivalence, only five such arrangements. Moreover, for reduced plane curves with quasi-homogeneous singularities, Schenck, Terao and Yoshinaga in \cite{STY}  present, among many results, an addition-deletion theorem that allows us to study their freeness. It is worth pointing out that arrangements of rational curves show up in the context of constructing ball-quotient surfaces or generally in uniformization of surfaces. Probably the most crucial results are those devoted to line arrangements in the complex projective plane which are reported in the classical textbook by Barthel, Hirzebruch and H\"ofer \cite{BHH87}. On the other hand, as it is explained by Kobayashi in \cite[Section 4.4]{Kob}, rational curves can be used to construct interesting counterexamples to uniformization of surfaces with $2c_{2}=c_{1}^{2}$. This small sample shows an undoubted importance of rational curve arrangements in the general theory of algebraic surfaces.

The main aim of the present note is to explore the world of conic arrangements in the complex projective plane with certain quasi-homogeneous singularities. In a very recent paper \cite{DJP}, the author with Dimca and Janasz observed that conic arrangements with nodes and tacnodes are never free. Here the point is to consider a wider class of conic arrangements that we are going to define right now.
\begin{definition}
\label{quass}
Let $\mathcal{C} = \{C_{1}, ..., C_{k}\} \subset \mathbb{P}^{2}_{\mathbb{C}}$ be an arrangement of $k\geq 2$ smooth conics. Assume that $\mathcal{C}$ has $n_{2}$ nodes, $t_{2}$ tacnodes, $n_{3}$ ordinary triple points, and $n_{4}$ ordinary quadruple points.
Then we call $\mathcal{C}$ as a $\mathcal{Q}$-conic arrangement.
\end{definition}
Of course, the letter $\mathcal{Q}$ refers to the fact that all prescribed singularities are \textit{quasi-homogeneous}. It is very hard to describe which topological types of singularities, in the world of conic arrangements, are quasi-homogeneous. It is worth pointing out here that for conic arrangements ordinary points of multiplicity greater or equal to $5$ are, in principle, not quasi-homogeneous. In order to explain this phenomenon, we are going to discuss two natural examples below. Before that, let us recall the following crucial definitions that we will used through the note.
\begin{definition}

Let $p$ be an isolated singularity of a polynomial $f\in \mathbb{C}[x,y]$. Since we can change the local coordinates, let $p=(0,0)$.
The number 
$$\mu_{p}=\dim_\mathbb{C}\left(\mathbb{C}[x,y] /\bigg\langle \frac{\partial f}{\partial x},\frac{\partial f}{\partial y} \bigg\rangle\right)$$
is called the Milnor number of $f$ at $p$.

The number
$$\tau_{p}=\dim_\mathbb{C}\left(\mathbb{C}[ x,y] /\bigg\langle f,\frac{\partial f}{\partial x},\frac{\partial f}{\partial y}\bigg\rangle \right)$$
is called the Tjurina number of $f$ at $p$.
\end{definition}
For a projective situation, with a point $p\in \mathbb{P}^{2}_{\mathbb{C}}$ and a homogeneous polynomial $f\in \mathbb{C}[x,y,z]$, we take local affine coordinates such that $p=(0,0,1)$ and then the dehomogenization of $f$.

\begin{definition}
A singularity is called quasi-homogeneous if and only if there exists a holomorphic change of variables so that the defining equation becomes weighted homogeneous.
\end{definition}
 Recall that $f(x,y) = \sum_{i,j}c_{i,j}x^{i}y^{j}$ is weighted homogeneous if there exist rational numbers $\alpha, \beta$ such that $\sum_{i,j} c_{i,j}x^{i\cdot \alpha} y^{j \cdot \beta}$ is homogeneous. It was proved by Reiffen in \cite{Reiffen} that if $f(x,y)$ is a convergent power series with an isolated singularity at the origin, then $f(x,y)$ is in the ideal generated by the partial derivatives if and only if $f$ is quasi-homogeneous. It means that in the quasi-homogeneous case one has $\tau_{p} = \mu_{p}$.
\begin{example}
Consider the arrangement $\mathcal{C} =\{C_{1}, ..., C_{5}\} \subset\mathbb{P}^{2}_{\mathbb{C}}$ given by
\begin{equation*}
\begin{array}{lr}
C_{1} : &(x-3z)^2 + (y-4z)^2 - 25z^2 = 0, \\
C_{2} : &(x-4z)^2 + (y-3z)^2 - 25z^2= 0, \\
C_{3} : &(x+3z)^2 + (y-4z)^2 - 25z^2 = 0, \\
C_{4} : &(x+4z)^2 + (y-3z)^2 - 25z^2 = 0, \\
C_{5} : &(x-5z)^2 + y^2 - 25z^2 = 0. \\
\end{array}
\end{equation*}
By \cite[Example 4.2]{SchenckToh} we know that the intersection point $p=(0:0:1)$, which is ordinary and has multiplicity $5$, is not quasi-homogeneous since $\tau_{p}=15$ but $\mu_{p}=16$.
\end{example}
\begin{example}
If we have an arrangement of $m$ reduced elements of the pencil of plane curves $f_{t} = g_{1} + tg_{2}$, where ${\rm deg} \, g_{1}={\rm deg} \, g_{2} = k$ such that the base locus $B$ consists of $k^2$ points, then any member of $B$ is an ordinary quasi-homogeneous singular point - see \cite{STY} for necessary details.
\end{example}
The discussion presented above allows us to justify the mentioned slogan that, in principle, ordinary singularities of multiplicity greater of equal to $5$ are not quasi-homogeneous.

Now we would like to pass to $\mathcal{Q}$-arrangements with singularities prescribed in Definition \ref{quass}. Our decision to choose such singularities is based on our experience with arrangements of lines and conics - we expect that such quasi-homogeneous singularities might be meaningful for the freeness problem.

First of all, we can observe that for $\mathcal{Q}$-conic arrangements $\mathcal{C}$ the following combinatorial count holds.
\begin{equation}
    4\cdot \binom{k}{2} = n_{2} + 2t_{2} + 3n_{3} + 6n_{4}. 
\end{equation}
\begin{proof}
The left-hand side is equal to the number of pairwise intersections of $k$ conics in $\mathcal{C}$. The right-hand side comes from the count of the intersection indices. Indeed, each node has the intersection index equal to $1$, tacnodes has the intersection index equal to $2$, and for ordinary triple and quadruple point we have $3$ and $6$, respectively.
\end{proof}
 Based on the above discussion regarding the freeness of reduced curves, it is very natural to wonder whether we can construct examples of free $\mathcal{Q}$-conic arrangements. In the class of conic arrangements, we are aware of one highly non-trivial arrangement called Chilean. It consists of $12$ conics and it has $9$ points of multiplicity $8$ and $12$ double points \cite{Dolg}. As it turned out, the Chilean arrangement is free \cite{PokSz}. However, the mentioned eightfold intersection points, even if these are ordinary singularities, are not quasi-homogeneous. 
 
 Turning back to our $\mathcal{Q}$-conic arrangements, we are able, to our surprise, prove the following result.
\begin{theoremA}
There does not exist any $\mathcal{Q}$-conic arrangement which is free.
\end{theoremA}
The above result stands in a strong opposition with respect to what we know about line arrangements, i.e, there are free line arrangements with double, triple, and quadruple points.

Looking at $\mathcal{Q}$-conic arrangements, it is natural to ask whether there are some additional constraints, different than the combinatorial count, on the weak combinatorics. Let us recall that by  the weak combinatorics of a given $\mathcal{Q}$-conic arrangement we understand the vector $(k;n_{2},t_{2},n_{3},n_{4})$, i.e., we focus only on the numerical data associated with $\mathcal{C}$, and we ignore the position of those singularities as the intersection points.
It turns out that we can show the following result which is in the spirit of results from \cite{Pokora2}.
\begin{theoremB}
Let $\mathcal{C}$ be a $\mathcal{Q}$-conic arrangement with $k\geq 3$. Then one has
$$8k + n_{2} + \frac{3}{4}n_{3} \geq \frac{5}{2}t_{2}.$$
\end{theoremB}
The note is organized as follows. In Section $2$, we recall all necessary definitions regarding Theorem $A$ and we present its proof. In Section $3$, we recall basics regarding an orbifold Miyaoka-Yau inequality and we present our proof of Theorem $B$. 

Through the paper we work over the complex numbers $\mathbb{C}$.

\section{Non-freeness of $\mathcal{Q}$-conic arrangements}
We start with the notion of free curves. Let $C$ be a reduced curve $\mathbb{P}^{2}_{\mathbb{C}}$ of degree $d$ defined by $f \in S :=\mathbb{C}[x,y,z]$. We denote by $J_{f}$ the Jacobian ideal generated by $\partial_{x}f, \, \partial_{y}f, \, \partial_{z}f$. We define by $r={\rm mdr}(f)$ the minimal degree of relations among the partial derivatives, i.e., the minimal degree $r$ of a triple $(a,b,c) \in S_{r}^{3}$ such that 
$$a\cdot \partial_{x} f + b\cdot \partial_{y}f + c\cdot \partial_{z}f = 0.$$
Consider now the graded $S$-module $N(f) = I_{f} / J_{f}$, where $I_{f}$ is the saturation of $J_{f}$ with respect to the irrelevant ideal $\mathfrak{m}=\langle x,y,z\rangle$. 
\begin{definition}
A reduced plane curve $C$ is \emph{free} if $N(f) = 0$.
\end{definition}
In order to study the freeness of reduced plane curves, we will use the following result due to du Plessis and Wall \cite{duP}.

\begin{theorem}[du Plessis - Wall]
\label{duW}
Let $C \subset \mathbb{P}^{2}_{\mathbb{C}}$ be a reduced plane curve of degree $d$ and let $f=0$ be its defining equation. Denote by $r: = {\rm mdr}(f)$ and assume that $r \leq (d-1)/2$. Then $C$ is free if and only if
\begin{equation}
\label{Milnor1}
r^2 - r(d-1) + (d-1)^2 = \tau(C),
\end{equation}
where $\tau(C) := \sum_{p \in {\rm Sing}(C)}\tau_{p}$ is the total Tjurina number of $C$.
\end{theorem}

Observe that if $C$ is a $\mathcal{Q}$-conic arrangement, then its global Tjurina number has the following form
$$\tau(C) = n_{2} + 3t_{2} + 4n_{3}  + 9n_{4}.$$
This formula follows from the fact that all these singularities are quasi-homogeneous and thus we have $$\tau(C) = \sum_{p \in {\rm Sing}(C)}\mu_{p}.$$

Now we are ready to prove Theorem A.
\begin{proof}
By Theorem \ref{duW}, a $\mathcal{Q}$-conic arrangement $\mathcal{C}$ with $k\geq 2$ is free if and only if
$$r^{2}-r(2k-1)+(2k-1)^{2} = n_{2}+3t_{2} +4n_{3} + 9n_{4}.$$
Using the combinatorial count, we get
$$r^{2} - r(2k-1)+(2k-1)^{2} = 4\cdot \binom{k}{2} + t_{2} + n_{3}  + 3n_{4}= 2(k^{2}-k) + t_{2} + n_{3}  + 3n_{4}.$$
This gives us 
\begin{equation}
\label{nf}
r^{2}-r(2k-1) + 2k^{2}-2k+1 -(t_{2} + n_{3}   + 3n_{4})=0.
\end{equation}
Now equation (\ref{nf}) has roots if
$$\triangle_{r} = (2k-1)^{2} - 4(2k^{2}-2k+1)+4(t_{2}+n_{3}  + 3n_{4})\geq 0.$$
This gives us the following inequality
\begin{equation}
    t_{2}+n_{3} + 3n_{4} \geq k^{2}-k+ \frac{3}{4}.
\end{equation}
Let us come back to the combinatorial count. Observe that
$$4\cdot \binom{k}{2} = 2k^{2} - 2k = n_{2} + 2t_{2} + 3n_{3} + 6n_{4} =$$ 
$$ n_{2} + n_{3} +  2(t_{2}+n_{3}  +3n_{4}) \geq n_{2} + n_{3} + 2k^{2} - 2k + \frac{3}{2},$$
and we arrive at $n_{2} + n_{3} \leq -\frac{3}{2}$ for every $k \in \mathbb{Z}_{\geq 2}$. This stands in a contradiction with the fact that $n_{2},n_{3}$ are non-negative, which completes the proof.

\end{proof}
\section{Combinatorial constraints on $\mathcal{Q}$-conic arrangements}
In order to prove Theorem B, we will follow Langer's variation on the Miyaoka-Yau inequality \cite{Langer} which uses the notion of local orbifold Euler numbers $e_{orb}$ of singular points. We are working with \emph{log pairs} $(X,D)$, where $X$ is a complex normal projective surface and $D$ is a boundary divisor which means that it is an effective $\mathbb{Q}$-divisor whose coefficients are less than or equal to one and such that $K_{X}+D$ is $\mathbb{Q}$-Cartier. Recall also that a log pair $(X,D)$ is effective if $K_{X}+D$ is effective.
\begin{definition}
Let $(X,D)$ be a log pair and let $f: Y \rightarrow X$ be a proper birational morphism from a normal surface $Y$. Write 
$$K_{Y}+D_{Y} = f^{*}(K_{X}+D)$$
with $f_{*}D_{Y} = D$. If the coefficients of $D_{Y}$ are less than or equal to one for every $f: Y \rightarrow X$, then $(X,D)$ is called a log canonical surface.
\end{definition}
In the light of the above definition, from now on we consider log pairs $(\mathbb{P}^{2}_{\mathbb{C}}, \alpha D)$, where $D$ is a boundary divisor consisting of $k$ smooth conics admitting only nodes, tacnodes, ordinary triple and quadruple points, and $\alpha$ is a suitable chosen positive rational number less than or equal to one.

Now we need to recall certain results regarding local orbifold numbers $e_{orb}$ that appear in the context of $\mathcal{Q}$-conic arrangements. We have to point out here that the definition of local orbifold Euler numbers is technical and it requires a severe preparation. Due to this reason, we aim to provide only some of their useful properties and for all necessary details we refer directly to \cite{Langer}. First of all, let us provide two comments that will shed some lights on these numbers, namely:
\begin{itemize}
\item local orbifold Euler numbers are analytic in their nature,
\item if $(\mathbb{C}^{2},D)$ is log canonical at $0$ and ${\rm mult}_{0}(D)$ denotes the multiplicity of $D$ in $0$ ( i.e., this is a sum of multiplicities of irreducible components $D_{i}$ counted with appropriate multiplicities), then
$$e_{orb}(0, \mathbb{C}^{2},D) \leq (1-{\rm mult}_{0}(D)/2)^{2},$$
which means that $e_{orb}(x,X,D)\leq 1$ for any log canonical pair $(X,D)$.
\end{itemize}
Going into details for $\mathcal{Q}$-arrangements, the local orbifold Euler numbers are the following:
\begin{itemize}
    \item if $q$ is a node, then $e_{orb}(q,\mathbb{P}^{2}_{\mathbb{C}}, \alpha D)=(1-\alpha)^2$ with $0 \leq \alpha \leq 1$;
    \item if $q$ is a tacnode, then $e_{orb}(q,\mathbb{P}^{2}_{\mathbb{C}}, \alpha D) =\frac{(3-4\alpha)^{2}}{8}$ with $\frac{1}{4} <  \alpha \leq \frac{3}{4}$;
    \item if $q$ is an ordinary triple or quadruple point, then $e_{orb}(q,\mathbb{P}^{2}_{\mathbb{C}}, \alpha D) \leq \bigg(1 - \frac{m\alpha}{2}\bigg)^2$ with $m \in \{3,4\}$ and   $0 \leq \alpha \leq \frac{2}{m}$.
\end{itemize}
Using the above discussion we can present our proof of Theorem B.
\begin{proof}
Let $\mathcal{C}$ be a $\mathcal{Q}$-conic arrangement and denote by $D:= C_{1} + ... + C_{k}$ the associated divisor to $\mathcal{C}$. In the first step, we determine $\alpha$ such that $K_{\mathbb{P}^{2}_{\mathbb{C}}} + \alpha C$ is effective and log canonical. For being effective, one needs to satisfy the condition that $-3 + 2k \alpha \geq 0$ and it implies that $\alpha \geq 3/2k$. Next, our pair is log-canonical if $\alpha \leq {\rm min} \{1, 3/4, 2/3, 1/2\}$, so $\alpha \leq 1/2$. Moreover, in order to apply the above formulae for local orbifold Euler numbers, we have to assume that $\alpha > 1/4$. Due to these reasons, we firstly arrive at $\alpha \in [3/2k, 1/2]$, which implies that $k\geq 3$, and by the constraint $\alpha > 1/4$, we select $\alpha=\frac{1}{2}$.  From now on we focus on the pair $\bigg(\mathbb{P}^{2}_{\mathbb{C}},\frac{1}{2}D\bigg)$, and we are going to use the inequality from \cite{Langer} adapted to our situation, namely
\begin{equation}
\label{logMY}
\sum_{p \in {\rm Sing}(\mathcal{C})}  3\bigg( \frac{1}{2}\bigg(\mu_{p} - 1\bigg) + 1 - e_{orb}\bigg(p,\mathbb{P}^{2}_{\mathbb{C}}, \frac{1}{2} D\bigg) \bigg) \leq 5k^{2} - 3k,
\end{equation}
where $\mu_{p}$ is the Milnor number of a singular point $p \in {\rm Sing}(\mathcal{C})$.

We can observe, after some straightforward calculations, that the left-hand side of \eqref{logMY} is bounded from below by
$$3n_{2}\cdot \bigg(1-\frac{1}{4}\bigg) + 3t_{2}\cdot \bigg(1+1-\frac{1}{8}\bigg) + 3n_{3}\cdot \bigg(\frac{3}{2}+1 - \frac{1}{16}\bigg) + 3n_{4}\cdot(4+1) = $$ 
$$\frac{9}{4}n_{2}+\frac{45}{8}t_{2}+\frac{117}{16}n_{3}+15n_{4}.$$
Since $$2k^{2} = 2k + n_{2} + 2t_{2} +3n_{3}+6n_{4},$$ we obtain
\begin{equation}
\label{ste}
\frac{9}{4}n_{2}+\frac{45}{8}t_{2}+\frac{117}{16}n_{3}+15n_{4} \leq \frac{5}{2}\bigg(2k + n_{2} + 2t_{2} +3n_{3}+6n_{4}\bigg)-3k.
\end{equation}
Now we multiply \eqref{ste} by $16$, and we arrive at
$$36n_{2}+90t_{2}+117n_{3}+240n_{4} \leq 32k + 40n_{2}+80t_{2}+120n_{3}+240n_{4}.$$
From this, we get
$$8k + n_{2} + \frac{3}{4}n_{3}\geq \frac{5}{2}t_{2},$$
which completes the proof.
\end{proof}

\begin{remark}
One may wonder why in our inequality there is no information about ordinary quadruple points, but this is a typical phenomenon when one applies variants of  Miyaoka-Yau inequality -- you can confront it with the celebrated Hirzebruch's inequality for line arrangements.
\end{remark}
\begin{remark}
Let us now focus on the case when our $\mathcal{Q}$-conic arrangement admits only nodes and tacnodes. Then we have $2k^{2}-2k = n_{2}+2t_{2}$ and our inequality above gives
$$8k + n_{2} \geq \frac{5}{2}t_{2}.$$
Combining the facts above
$$8k + 2k^{2}-2k -2t_{2} \geq \frac{5}{2}t_{2},$$
we arrive at
$$t_{2} \leq \frac{4}{9}k^{2} + \frac{4}{3}k,$$
obtaining the same upper-bound on the number of tacnodes for conic arrangements as in Miyaoka's paper \cite{Miyaoka}.
\end{remark}

\section*{Acknowledgments}
The author would like to thank Alex Dimca for many discussions regarding conic arrangements in the plane and for encouraging to write down the present note. Moreover, he would like to express his gratitude to an anonymous referee for very helpful comments that allowed to improve the exposition of the note.

The author was partially supported by the National Science Center (Poland) Sonata Grant Nr \textbf{2018/31/D/ST1/00177}.

\vskip 0.5 cm

Piotr Pokora,
Department of Mathematics,
Pedagogical University of Krakow,
Podchor\c a\.zych 2,
PL-30-084 Krak\'ow, Poland. \\
\nopagebreak
\textit{E-mail address:} \texttt{piotr.pokora@up.krakow.pl}
\bigskip
\end{document}